\documentclass[leqno]{amsart}
\usepackage{amsmath}
\usepackage{amssymb}
\usepackage{amsthm}
\usepackage{enumerate}
\usepackage[mathscr]{eucal}
\theoremstyle{plain}
\usepackage{tikz}
\newtheorem{theorem}{Theorem}[section]
\newtheorem{lemma}[theorem]{Lemma}
\newtheorem{prop}[theorem]{Proposition}
\theoremstyle{definition}
\newtheorem{definition}[theorem]{Definition}
\newtheorem{remark}[theorem]{Remark}

\newtheorem{example}[theorem]{Example}
\newtheorem*{acknowledgement}{\textnormal{\textbf{Acknowledgements}}}

\newtheorem{cor}[theorem]{Corollary}
\theoremstyle{remark}

\begin{document}

\title{On Best approximations to compact operators}

\author{Debmalya Sain}

\address{(Sain)~Department of Mathematics, Indian Institute of Science, Bengaluru 560012, Karnataka, India}
\email{saindebmalya@gmail.com}


\thanks{}

\subjclass[2010]{Primary 46B28,  Secondary 46B20}
\keywords{compact operators; Birkhoff-James orthogonality; best approximations; operator norm attainment}



\date{}

\begin{abstract}
 We study best approximations to compact operators between Banach spaces and Hilbert spaces, from the point of view of Birkhoff-James orthogonality and semi-inner-products. As an application of the present study, some distance formulae are presented in the space of compact operators. The special case of bounded linear functionals as compact operators is treated separately and some applications to best approximations in reflexive, strictly convex and smooth Banach spaces are discussed. An explicit example is presented in $ \ell_p^{n} $ spaces, where $ 1 < p < \infty, $ to illustrate the applicability of the methods developed in this article. A comparative analysis of the results presented in this article with the well-known classical duality principle in approximation theory is conducted to demonstrate the advantage in the former case, from a computational point of view.
\end{abstract}

\maketitle 

\section{Introduction}
The purpose of this article is to study best approximations in the space of compact operators between Banach spaces and Hilbert spaces, and to present some distance formulae in certain special cases.  Let us first establish the notations and the terminologies to be used throughout the article.\\
The letters $\mathbb{X}, \mathbb{Y}, \mathbb{Z}$ stand for Banach spaces and the letter $ \mathbb{H} $ is used to denote a Hilbert space. The symbol $ \langle~,~\rangle $ is used to denote the inner product on $ \mathbb{H}. $  Let $ \theta $ denote the zero vector of any vector space, other than the scalar field. We work with both real and complex Hilbert spaces and \emph{only} real Banach spaces. Let $ B_{\mathbb{X}}= \{ x \in \mathbb{X}~:\|x\| \leq 1 \} $ and $  S_{\mathbb{X}}= \{ x \in \mathbb{X}~:\|x\| = 1 \} $ be the unit ball and the unit sphere of $ \mathbb{X}, $ respectively. $ \mathbb{X} $ is said to be strictly convex if every element of $ S_{\mathbb{X}} $ is an extreme point of $ B_{\mathbb{X}}. $ Let $K(\mathbb{X}, \mathbb{Y})$ denote the Banach space of all compact operator from $\mathbb{X}$ to $ \mathbb{Y}, $ endowed with the usual operator norm and let $\mathbb X^*$ denote the dual space of $\mathbb{X}.$ For the sake of brevity, we write $ K(\mathbb{X}, \mathbb{Y}) = K(\mathbb{X}), $ whenever $ \mathbb{X}=\mathbb{Y}. $ Given $ T \in K(\mathbb{X}, \mathbb{Y}), $ we use the notations $ \mathcal{R}(T) $ and $ \mathcal{N}(T) $ to denote the range of $ T $ and the kernel of $ T, $ respectively. The study of best approximation(s) to a given element out of a given subspace is a classical area of research in functional analysis. Let us recall the following basic definition in this context:

\begin{definition}
Given an element $ x \in \mathbb{X} $ and a subspace $ \mathbb{Y} $ of $ \mathbb{X}, $ let $ \textit{dist}\{x, \mathbb{Y}\} = \inf\{ \| x - z \| : z \in \mathbb{Y} \} $ denote the distance between $ x $ and $ \mathbb{Y}. $ An element $ y \in \mathbb{Y} $ is said to be a best approximation to $ x $ out of $ \mathbb{Y} $ if $ \| x-y \| = \min\{ \| x-z \| : z \in \mathbb{Y} \}.$
\end{definition}

It is well-known that in general neither the existence nor the uniqueness of best approximation is guaranteed. However, best approximation(s) always exist for finite-dimensional subspaces, and more generally for compact subsets of infinite-dimensional subspaces. Moreover, best approximation is unique in a strictly convex Banach space, provided it exists. Clearly, the above definition makes sense in the space of operators between Banach (Hilbert) spaces, and it is worth studying only when $ x \notin \mathbb{Y}. $ For the study of best approximations of linear operators in various contexts and under additional assumptions, one may consult \cite{A, AR, W} and the references therein. The study of best approximations in Banach spaces is intimately connected to the concepts of Birkhoff-James orthogonality and semi-inner-products. Given $ x, y \in \mathbb{X}, $ we say that $ x $ is Birkhoff-James orthogonal to $ y, $ written as $ x \perp_B y, $ if $ \| x+\lambda y \| \geq \| x \| $ for all scalars $ \lambda. $ It is easy to observe that Birkhoff-James orthogonality is homogeneous, i.e., $ x \perp_B y $ implies that $ \alpha x \perp_B \beta y $ for all scalars $ \alpha, \beta. $ Moreover, we note that in a Hilbert space, the Birkhoff-James orthogonality relation $ \perp_B $ coincides with the usual orthogonality relation $ \perp $ induced by the underlying inner product $ \langle~,~\rangle. $ Following \cite{S}, we say that $ y \in x^+ (y \in x^-) $ if $ \| x+\lambda y \| \geq \| x \| $  for all $ \lambda \geq 0 (\lambda \leq 0). $ We use the notation $ x^{\perp} = \{ y \in \mathbb{X} : x \perp_B y \} $ to denote the Birkhoff-James orthogonality set of the vector $ x. $ We refer the readers to the pioneering articles \cite{B, J, Ja} for the basic applications of Birkhoff-James orthogonality in understanding the geometry of Banach spaces, and to \cite{BS, S, SP, SPM, SPMa} for some of the more recent works in Banach spaces, involving the said notion of orthogonality. It is easy to observe that $ y \in \mathbb{Y} $ is a best approximation to $ x $ out of $ \mathbb{Y} $ if and only if $ (x-y) \perp_B \mathbb{Y}, $ i.e., $ (x-y) \perp_B z $ for all $ z \in \mathbb{Y}. $ We next mention the concept of semi-inner-products in Banach spaces, which is integral to the theme of this article.

\begin{definition}
Let $ \mathbb{X} $ be a real Banach space. A function $ [ ~,~ ] : \mathbb{X} \times \mathbb{X} \longrightarrow \mathbb{R} $ is a semi-inner-product (s.i.p.) if for any $ \alpha,~\beta \in \mathbb{R} $ and for any $ x,~y,~z \in \mathbb{X}, $ it satisfies the following:\\
$ (a) $ $ [\alpha x + \beta y, z] = \alpha [x,z] + \beta [y,z], $\\
$ (b) $ $ [x,x] > 0, $ whenever $ x \neq 0, $\\
$ (c) $ $ |[x,y]|^{2} \leq [x,x] [y,y], $\\
$ (d) $ $ [x,\alpha y] = \alpha [x,y]. $
\end{definition} 

It was proved in \cite{G} (see also \cite{Sa} for a rigorous proof of the same) that every Banach space $ (\mathbb{X},\|.\|) $ can be represented as an s.i.p. space $ (\mathbb{X},[~,~]) $ such that for all $ x \in \mathbb{X}, $ it holds that $ [x,x] = \| x \|^{2}. $ Whenever we speak of an s.i.p. $ [~,~] $ in the context of a Banach space $ \mathbb{X} $, we implicitly assume that  $ [~,~] $ is compatible with the norm, i.e., for all $ x \in \mathbb{X}, $ we have, $ [x,x] = \| x \|^{2}. $ In general, there can be many compatible s.i.p. corresponding to a given norm. As observed by Lumer in \cite{L}, there exists a unique s.i.p. on a normed space if and only if the space is smooth. We recall that $ \mathbb{X} $ is said to be smooth if there exists a unique supporting hyperplane to $ B_{\mathbb{X}} $ at each point of $ S_{\mathbb{X}}. $\\ 

Bhatia and $ \breve{S} $emrl studied Birkhoff-James orthogonality of matrices (viewed as operators on a finite-dimensional Hilbert space) and obtained some useful distance formulae in their seminal article \cite{BS}. Indeed, Theorem $ 1.1 $ of \cite{BS}, also known as the Bhatia-$\breve{S}$emrl Theorem, gives a complete characterization of the Birkhoff-James orthogonality of matrices. For an analogous study of orthogonality of operators between real Banach spaces, we refer the readers to \cite{S, SP, SPM}. For the study of orthogonality and best approximations in the space of matrices (a special case of compact operators on a Hilbert space), one should see \cite{Ga, LR}. The notion of the norm attainment set of an operator plays a central role in the study of orthogonality of operators. Given $ T \in K(\mathbb{X}, \mathbb{Y}), $ let $ M_T=\{ x \in S_{\mathbb{X}} : \| Tx \|=\|T\| \} $ denote the norm attainment set of the operator $ T. $ Observe that whenever $ \mathbb{X} $ is reflexive, it follows that $ M_T \neq \emptyset. $\\

Our aim in the present article is to further build upon the ideas presented in the above mentioned works, in the space of compact operators. We obtain a complete characterization of best approximations to a given compact operator out of a given finite-dimensional subspace, separately for Hilbert spaces and reflexive Banach spaces. In particular, distance formula for a compact operator and a one-dimensional subspace are presented in both the cases. As the most important part of this article, we treat the special case of functionals as compact operators and present an efficient algorithm to study the best approximation problems in reflexive smooth and strictly convex Banach spaces. Explicit examples are presented in the setting of $ \ell_p^n $ spaces to illustrate the applicability of the methods developed here. In the short final section of this article, we make a comparative analysis of our results with the classical duality principle in approximation theory. Indeed, we show that from a purely computational point of view, the use of orthogonality can evidently strengthen the well-known duality principle. 

\section{Best approximations to compact operators}
We begin with the observation that Bhatia-$\breve{S}$emrl type theorems for compact operators on a reflexive Banach space immediately give a complete characterization of best approximations out of a one-dimensional subspace.

\begin{prop}\label{prop:one-dim-Banach}
Let $ \mathbb{X} $ be a reflexive Banach space and let $ \mathbb{Y} $ be any Banach space. Let $ T, A \in K(\mathbb{X}, \mathbb{Y}) $ be linearly independent and let $ \lambda_0 \in \mathbb{R}. $ Then the following are equivalent:
\item[(i)] $ (T-\lambda_0 A) \perp_B A, $
\item[(ii)] there exist $ x, y \in M_{T-\lambda_0 A} $ and s.i.p. $ [~,~]_1, [~,~]_2 $ on $ \mathbb{Y} $ such that $ [Ax, Tx-\lambda_0 Ax]_1
\geq 0 $ and $ [Ay, Ty-\lambda_0 Ay]_2\leq 0,  $
\item[(iii)] $ \lambda_0 A $ is a best approximation to $ T $ out of $ \textit{span}\{A\}. $
\end{prop}

\begin{proof}
The equivalence of (i) and (iii) follows from the definitions of Birkhoff-James orthogonality and best approximations. Let us first prove that (i) implies (ii). It follows from Theorem $ 2.1 $ of \cite{SPM} that there exist $ x, y \in M_{T-\lambda_0 A} $ such that $ Ax \in ((T-\lambda_0 A)x)^{+} $ and $ Ay \in ((T-\lambda_0 A)y)^{-}. $ Applying Theorem $ 2.4 $ of \cite{SPMa}, we deduce that (ii) holds true. In similar spirit, applying the converses of these two theorems, we obtain that (ii) implies (i).
\end{proof}

If $ M_{T-\lambda_0 A} $ is of a particularly nice form then we have a refinement of the above observation. The proof is omitted as it follows directly from Theorem $ 2.2 $ of \cite{SPMa} and the observation that $ (T-\lambda_0 A)x \perp_B Ax $ if and only if there exists an s.i.p. $ [~,~] $ on $ \mathbb{Y} $ such that $ [Ax, Tx-\lambda_0 Ax]=0. $

\begin{prop}\label{prop:one-dim-Banach*}
Let $ \mathbb{X} $ be a reflexive Banach space and let $ \mathbb{Y} $ be any Banach space. Let $ T, A \in K(\mathbb{X}, \mathbb{Y}) $ be linearly independent and let $ \lambda_0 \in \mathbb{R}. $ Also assume that $ M_{T-\lambda_0 A} = \pm D, $ where $ D $ is a connected subset of $ S_{\mathbb{X}}. $ Then the following are equivalent:
\item[(i)] $ (T-\lambda_0 A) \perp_B A, $
\item[(ii)] there exists $ x \in M_{T-\lambda_0 A} $ and an s.i.p. $ [~,~] $ on $ \mathbb{Y} $ such that $ [Ax, Tx-\lambda_0 Ax] = 0, $
\item[(iii)] $ \lambda_0 A $ is a best approximation to $ T $ out of $ \textit{span}\{A\}. $
\end{prop}

In case of compact operators on a Hilbert space, we have yet another refinement of both the above results. 

\begin{prop}\label{prop:one-dim-Hilbert}
Let $ \mathbb{H} $ be a Hilbert space and let $ T, A \in K(\mathbb{H}) $ be linearly independent. Let $ \lambda_0 \in \mathbb{C}. $ Then the following are equivalent:
\item[(i)] $ (T-\lambda_0 A) \perp_B A, $
\item[(ii)] there exists $ x \in M_{T-\lambda_0 A} $ such that $ \langle Tx, Ax \rangle = \lambda_0 \| Ax \|^2, $
\item[(iii)] $ \lambda_0 A $ is a best approximation to $ T $ out of $ \textit{span}\{A\}. $
\end{prop}

\begin{proof}
To prove the equivalence of (i) and (ii), we first observe that the inner product $ \langle~,~\rangle $ is the only s.i.p. on $ \mathbb{H}. $ Applying Remark $ 3.1 $ of \cite{BS} (also see Theorem 2.2 of \cite{SPMa} and Theorem $ 2.2 $ of \cite{SP} for the real case), and the compactness of the operator $ T-\lambda_0 A, $ we deduce that $ (T-\lambda_0 A) \perp_B A $ if and only if there exists $ x \in M_{T-\lambda_0 A} $ such that $ \langle Ax, Tx-\lambda_0Ax \rangle =0. $ Now the desired conclusion follows from the conjugate symmetry and additivity properties of the inner product $ \langle~,~\rangle $.
\end{proof}

\begin{remark}
In view of the above results, a natural question arises regarding the uniqueness of best approximation to a compact operator out of a one-dimensional subspace. This can be answered in terms of a strengthening of Birkhoff-James orthogonality introduced in \cite{PSJ}. Given $ x, y \in \mathbb{X}, $ we say that $ x $ is strongly orthogonal to $ y $ in the sense of Birkhoff-James, written as $ x \perp_{SB} y, $ if $ \| x+\lambda y \| > \| x \| $ for all $ \lambda \neq 0. $ In each of the above propositions, it is easy to observe that $ \lambda_0 A $ is the unique best approximation to $ T $ out of $ \textit{span}\{A\} $ if and only if $ (T-\lambda_0 A) \perp_{SB} A. $
\end{remark}

In Theorem $ 2.9 $ of \cite{SPMa}, a distance formula has been presented for a compact operator and a one-dimensional subspace, under additional assumptions on norm attainment of certain operators. We would like to observe that the norm attainment condition can be relaxed, without any essential changes to the argument presented there. To this end, we first prove the following modification of Theorem $ 2.5 $ of \cite{SPMa}.

\begin{lemma}\label{lemma:norm evaluation}
Let $ \mathbb{X} $ be a reflexive Banach space and let $ \mathbb{Y} $ be any Banach space. Let $ T, A \in K(\mathbb{X}, \mathbb{Y}) $ be such that $ T \perp_B A $ and $ M_T = \pm D, $ where $ D $ is a connected subset of  $ S_{\mathbb{X}}. $ Then there exists an s.i.p. $ [~,~] $ on $ \mathbb{Y} $ such that
\begin{align*}
\| T \| & = \max\{|[Tx, y]| : x \in S_{\mathbb{X}}, y \in S_{\mathbb{Y}}, y \perp_B Ax\} \\
        & =  \max\{|[Tx, y]| : x \in S_{\mathbb{X}}, y \in S_{\mathbb{Y}},  [Ax, y] = 0\}.
\end{align*}
\end{lemma}

\begin{proof}
It follows from the defining properties of s.i.p. that given any s.i.p. $ [~,~] $ on $ \mathbb{Y}, $ the following holds true:
\[ \max\{|[Tx, y]| : x \in S_{\mathbb{X}}, y \in S_{\mathbb{Y}}, y \perp_B Ax\} \leq \|T\|. \]

On the other hand, Theorem $ 2.2 $ of \cite{SPMa} implies that there exists $ x_0 \in M_T $ such that $ Tx_0 \perp_B Ax_0. $ If $ Tx_0 = \theta $ then $ T $ is the zero operator and the result follows trivially. Let us assume that $ Tx_0 \neq \theta. $  Now, taking $ x=x_0 $ and $ y=\frac{Tx_0}{\|Tx_0\|}, $ we obtain by the homogeneity of Birkhoff-James orthogonality that $ x \in S_{\mathbb{X}}, y \in S_{\mathbb{Y}} $ and $ y \perp_B  Ax. $ Since 
\[ |[Tx, y]| = | [Tx_0, \frac{Tx_0}{\|Tx_0\|} ] | = \| Tx_0 \| = \| T \|, \]
the first equality follows. To prove the second equality, we only need to observe that $ y \perp_B Ax $ if and only if there exists an s.i.p. $ [~,~] $ on $ \mathbb{Y} $ such that $ [Ax, y] = 0. $ Indeed, applying Theorem $ 2.1 $ of \cite{J}, there exists $ f \in S_{\mathbb{Y}^*} $ such that $ f(y)=\|y\|=1 $ and $ f(Ax)=0. $ Following the argument presented by Lumer in Theorem $ 2 $ of \cite{L}, we can construct an s.i.p. $ [~,~] $ on $ \mathbb{Y} $ such that $ [Ax, y] = f(Ax) = 0. $ 
This establishes the lemma.
\end{proof}

\begin{remark}\label{remark:Hilbert}
We would like to note that the above lemma holds true for compact operators on a Hilbert space, real or complex, with the s.i.p. replaced by the underlying inner product $ \langle~,~\rangle. $ This can be verified easily by an easy application of Remark $ 3.1 $ of \cite{BS}.
\end{remark}

Applying the above lemma, we next obtain an improvement of Theorem $ 2.9 $ of \cite{SPMa}.

\begin{theorem}\label{theorem:distance formula}
Let $ \mathbb{X} $ be a reflexive Banach space and let $ \mathbb{Y} $ be any Banach space. Let $ T, A \in K(\mathbb{X}, \mathbb{Y}) $ be such that $ T \notin \textit{span}\{ A \}. $ Let $ \lambda_0 A $ be a best approximation to $ T $ out of $ \textit{span}\{A\}. $ Also assume that $ M_{T-\lambda_0 A} = \pm D, $ where $ D $ is a connected subset of $ S_{\mathbb{X}}. $ Then there exists an s.i.p. $ [~,~] $ on $ \mathbb{Y} $ such that 
$$\textit{dist}(T, \textit{span}\{ A \})=\max\{|[Tx, y]| : x \in S_{\mathbb{X}}, y \in S_{\mathbb{Y}}, [Ax, y] = 0\}.$$
\end{theorem}

\begin{proof}
As $ T-\lambda_0 A, A $ are compact, $ (T-\lambda_0 A) \perp_B A $ and $ M_{T-\lambda_0 A} = \pm D, $ where $ D $ is a connected subset of $ S_{\mathbb{X}}, $ it follows  from Lemma \ref{lemma:norm evaluation} that there exists an s.i.p. $ [~,~] $ on $ \mathbb{Y} $ such that $ \| T-\lambda_0 A \| = \max\{ |[(T-\lambda_0 A)x, y]| : x \in S_{\mathbb{X}}, y \in S_{\mathbb{Y}}, [Ax, y]=0 \}. $ Now, the desired result follows from the properties of s.i.p. and the fact that $ \| T-\lambda_0 A \| = \textit{dist}(T, \textit{span}\{ A \}). $ This establishes the theorem. 
\end{proof}

The above distance formula becomes more useful from a computational point of view when $ [~,~] $ is the unique s.i.p. on $ \mathbb{Y}, $ or, equivalently, when $ \mathbb{Y} $ is smooth. In particular, given any $ A \in K(\mathbb{X}, \mathbb{Y}), $ the distance of a compact operator $ T \in K(\mathbb{X}, \mathbb{Y}) $ from $ \textit{span}\{A\} $ can be obtained by computing the maximum of the quantities $ [Tx, y], $ where $ x \in S_{\mathbb{X}}, y \in S_{\mathbb{Y}}, $ and $ [Ax, y] = 0, $ provided the condition on the norm attainment set is satisfied. The interesting thing to observe in this context is that we do not require the smoothness condition on $ \mathbb{X}. $ This is illustrated in the following example:

\begin{example}
Let $ T, A \in K(\ell_{\infty}^2, \ell_{2}^2) $ be given by
$$ T(u, v) = (u+2v, 5u+5v), A(u, v) = (u, 0)~ \textit{for all}~ (u, v) \in \ell_{\infty}^2. $$

Let us compute $ \textit{dist}(T, \textit{span}\{ A \}), $ by applying Theorem \ref{theorem:distance formula}. Given any $ B \in K(\ell_{\infty}^2, \ell_{2}^2), $ we begin with the following two basic observations:\\

\noindent (i) $ \| B \|_{\infty, 2} = \max\{\|B(1, 1)\|_2, \|B(1, -1)\|_2\}, $\\
(ii) if $ \| B(1, 1) \|_2 \neq \| B(1, -1) \|_2 $ then $ M_B $ is necessarily of the form $ M_B = \pm D, $ where $ D $ is a connected subset of $ S_{\ell_{\infty}^2}. $ In fact, it is easy to show by using the convexity of norm that in this case $ D $ must be singleton.\\

Now, taking $ B = T-\lambda A, $ where $ \lambda \in \mathbb{R}, $ an easy computation reveals that $ \|  B(1, 1) \|_2 = \|  B(1, -1) \|_2 $ if and only if $ \lambda = \lambda_1 = 13.5.$ Let $ \lambda_0 A $ be a best approximation to $ T $ out of $ \textit{span}\{A\}. $ Since $ \| T-\lambda_1 A \|_{\infty, 2} = 14.5 > \| T \|_{\infty, 2} = \sqrt{109} $, hence $\lambda_1\neq \lambda_0,$ so for $B_0=T-\lambda_0A,$ we have $\|  B_0(1, 1) \|_2 \neq \|  B_0(1, -1) \|_2,$ and we conclude that the norm attainment condition on the operator $ T - \lambda_0 A $ in Theorem \ref{theorem:distance formula} is satisfied. Taking any $ x = (u, v) \in S_{\ell_{\infty}^2}, $ we note that $ Ax $ is a scalar multiple of $ (1, 0) \in S_{\ell_{2}^2}. $ Therefore, by the homogeneity of Birkhoff-James orthogonality, if $ y \in S_{\ell_{2}^2} $ is such that $ y \perp_B Ax, $ then $ y = \pm (0, 1). $ In particular, $ | [Tx, y] | = | [ (u+2v, 5u+5v), (0, 1) ] | = | 5u+5v |, $ where $ [~,~] = \langle~,~\rangle, $ the usual inner product on $ \mathbb{R}^2 .$ Maximizing over $ (u, v) \in S_{\ell_{\infty}^2}, $ we obtain by applying Theorem \ref{theorem:distance formula}, 
\[ \textit{dist}(T, \textit{span}\{ A \}) = 10.  \]
\end{example}

\begin{remark}
The above example shows that Theorem \ref{theorem:distance formula} is a proper improvement of Theorem $ 2.9 $ of \cite{SPMa}. Indeed, Theorem $ 2.9 $ of \cite{SPMa} is not applicable in case of the above example, as it can be easily checked that $ M_{T-\lambda_1 A} = \{ \pm (1, 1), \pm (1, -1) \}, $ which is not of the desired form.

\end{remark}

In case of compact operators on a Hilbert space, we have the following corollary to Theorem \ref{theorem:distance formula}, the proof of which is omitted as it is now obvious in light of Remark \ref{remark:Hilbert}. 

\begin{cor}
Let $ \mathbb{H} $ be a Hilbert space and let $ T, A \in K(\mathbb{H}) $ be such that $ T \notin \textit{span}\{ A \}. $ Then
$$\textit{dist}(T, \textit{span}\{ A \})=\max\{|\langle Tx, y \rangle | : x, y \in S_{\mathbb{H}}, y \perp Ax\}.$$
\end{cor}

\begin{remark}
When $ \mathbb{H} $ is finite-dimensional, the above corollary was proved in \cite{BS}, for the special case $ A=I, $ where $ I $ is the identity operator on $ \mathbb{H}. $ In this context, the significance of Theorem \ref{theorem:distance formula} is to illustrate that the concerned distance formula is valid for compact operators on a reflexive Banach space, with certain natural modifications. Moreover, the corresponding statement in case of Hilbert spaces follows directly from Theorem \ref{theorem:distance formula}.
\end{remark}

We next obtain a complete characterization of the best approximation to a compact operator out of an $ n- $dimensional subspace, where $ n \geq 2, $ under an additional condition on the norm attainment set.

\begin{theorem}\label{theorem:n-dim}
Let $ \mathbb{X} $ be a reflexive Banach space, let $ \mathbb{Y} $ be any Banach space and let $ n \geq 2. $ Let $ T, A_1, A_2, \ldots, A_n \in K(\mathbb{X}, \mathbb{Y}) $ be such that $ A_1, A_2, \ldots, A_n $ are linearly independent and $ T \notin \mathbb{Z}, $ where $ \mathbb{Z} = \textit{span}\{ A_1, A_2, \ldots, A_n \}. $ Also assume that $ M_{T-S} = \pm D_{S}, $ where $ D_{S} $ is a connected subset of $ S_{\mathbb{X}}, $ for every $ S \in \mathbb{Z}. $ Let $ \alpha_i \in \mathbb{R}, $ where $ 1 \leq i \leq n. $  Then $ \sum_{i=1}^{n} \alpha_i A_i $ is a best approximation to $ T $ out of $ \mathbb{Z} $ if and only if given any $ n $ scalars $ \beta_1, \beta_2, \ldots, \beta_n, $ there exist scalars $  \gamma_1, \gamma_2, \ldots, \gamma_n $ and s.i.p. $ [~,~]_1, [~,~]_2 $ on $ \mathbb{Y} $ such that the following holds true:
\begin{align*}
\| T - \sum_{i=1}^{n} \beta_i A_i \| & = \max\{ | [(T-\sum_{i=1}^{n} \beta_i A_i)x, y]_1 | : x \in S_{\mathbb{X}}, y \in S_{\mathbb{Y}}, y \perp_B \sum_{i=1}^{n} \gamma_i A_{i}x \} \\
 & \geq \max\{ | [(T-\sum_{i=1}^{n} \alpha_i A_i)x, y]_2 | : x \in S_{\mathbb{X}}, y \in S_{\mathbb{Y}}, y \perp_B \sum_{i=1}^{n} \gamma_i A_{i}x \} \\
 & = \| T - \sum_{i=1}^{n} \alpha_i A_i \|.
\end{align*}
\end{theorem}

\begin{proof}
As the sufficient part of the theorem follows trivially, we prove only the necessary part. An easy application of the Hahn-Banach Theorem shows that $ (T - \sum_{i=1}^{n} \beta_i A_i)^{\perp} $ contains a subspace of codimension one in $ K(\mathbb{X}, \mathbb{Y}). $ Since $ n \geq 2, $ it follows that $ (T - \sum_{i=1}^{n} \beta_i A_i)^{\perp} \bigcap \mathbb{Z} \neq \{ \theta \}. $ Let $ \sum_{i=1}^{n} \gamma_i A_{i} \in (T - \sum_{i=1}^{n} \beta_i A_i)^{\perp} \bigcap \mathbb{Z} $ be non-zero, where $ \gamma_i (1 \leq i \leq n) $ are scalars. Since $ (T - \sum_{i=1}^{n} \beta_i A_i) \perp_B \sum_{i=1}^{n} \gamma_i A_i, $ and $ M_{(T - \sum_{i=1}^{n} \beta_i A_i)} $ is of the desired form, applying Lemma \ref{lemma:norm evaluation}, we obtain that 
\[ \| T - \sum_{i=1}^{n} \beta_i A_i \|  = \max\{ | [(T-\sum_{i=1}^{n} \beta_i A_i)x, y]_1 | : x \in S_{\mathbb{X}}, y \in S_{\mathbb{Y}}, y \perp_B \sum_{i=1}^{n} \gamma_i A_{i}x \}, \] for some s.i.p. $ [~,~]_1 $ on $ \mathbb{Y}, $ proving the first equality. On the other hand, we also have that $ (T-\sum_{i=1}^{n} \alpha_i A_i) \perp_B \sum_{i=1}^{n} \gamma_i A_i \in \mathbb{Z}, $ as $ \sum_{i=1}^{n} \alpha_i A_i $ is a best approximation to $ T $ out of $ \mathbb{Z}. $ Therefore, applying similar arguments, we obtain that
\[ \| T - \sum_{i=1}^{n} \alpha_i A_i \|  = \max\{ | [(T-\sum_{i=1}^{n} \alpha_i A_i)x, y]_2 | : x \in S_{\mathbb{X}}, y \in S_{\mathbb{Y}}, y \perp_B \sum_{i=1}^{n} \gamma_i A_{i}x \}, \] for some s.i.p. $ [~,~]_2 $ on $ \mathbb{Y}, $ proving the last equality. The remaining inequality follows from the fact that $ \| T - \sum_{i=1}^{n} \beta_i A_i \| \geq \| T - \sum_{i=1}^{n} \alpha_i A_i \|. $ This establishes the theorem.
\end{proof}

In view of the above theorem, we note that in case $ \mathbb{Y} $ is smooth, $ [~,~]_1 = [~,~]_2. $ Regarding the uniqueness of best approximation, we make the following remark.

\begin{remark}
$ \sum_{i=1}^{n} \alpha_i A_i $ is the unique best approximation to $ T $ out of $ \mathbb{Z} = \textit{span}\{ A_1, A_2, \ldots, A_n \}, $ if and only if the inequality in Theorem \ref{theorem:n-dim} is strict. This can be verified easily by following the proof of the said theorem. 
\end{remark}

In case of compact operators on a Hilbert space, we have the following corollary to Theorem \ref{theorem:n-dim}. The proof is an easy adaptation of the proof of Theorem 2.12 since there is a unique s.i.p. in the Hilbert space.

\begin{cor}\label{corollary:n-dim-Hilbert}
Let $ \mathbb{H} $ be a Hilbert space and let $ T, A_1, A_2, \ldots, A_n \in K(\mathbb{H}) $ be such that $ A_1, A_2, \ldots, A_n $ are linear independent and $ T \notin \textit{span}\{ A_1, A_2, \ldots, A_n \}. $ Let $ \alpha_i \in \mathbb{C}, $ where $ 1 \leq i \leq n. $  Then $ \sum_{i=1}^{n} \alpha_i A_i $ is a best approximation to $ T $ out of $ \textit{span}\{ A_1, A_2, \ldots, A_n \} $ if and only if given any $ n $ scalars $ \beta_1, \beta_2, \ldots, \beta_n, $ there exist scalars $  \gamma_1, \gamma_2, \ldots, \gamma_n $ such that the following holds true:
\begin{align*}
\| T - \sum_{i=1}^{n} \beta_i A_i \| & = \max\{ | \langle (T-\sum_{i=1}^{n} \beta_i A_i)x, y \rangle | : x, y \in S_{\mathbb{H}}, y \perp \sum_{i=1}^{n} \gamma_i A_{i}x \} \\
 & \geq \max\{ | \langle(T-\sum_{i=1}^{n} \alpha_i A_i)x, y \rangle | : x, y \in S_{\mathbb{H}}, y \perp \sum_{i=1}^{n} \gamma_i A_{i}x \} \\
 & = \| T - \sum_{i=1}^{n} \alpha_i A_i \|.
\end{align*}
Moreover, the best approximation is unique if and only if the above inequality is strict.
\end{cor}

\section{Best approximations in reflexive spaces}

In light of the results obtained so far in this article, it is evident that for best approximations to compact operators, norm attainment set plays a central role. By virtue of Remark $ 1 $ in \cite{SP}, the norm attainment set of a compact operator on a Hilbert space $ \mathbb{H} $ is necessarily the unit sphere of some subspace of $ \mathbb{H}. $ Since such a nicety is no longer present in case of compact operators between Banach spaces, we require additional conditions in that case. However, for the class of bounded linear functionals on a reflexive Banach space, we do have an additional advantage as mentioned in the following proposition. The proof is easy, and is therefore omitted.

\begin{prop}\label{prop:functional}
Let $ \mathbb{X} $ be a reflexive Banach space and let $ f \in X^{*}. $ Then $ M_f $ is of the form $ \pm D, $ where $ D $ is a connected subset of $ S_{\mathbb{X}}$ ($D$ is a face of $B_\mathbb{X}$).
\end{prop}

Let us now present an analogous result to Proposition \ref{prop:one-dim-Banach} and Theorem \ref{theorem:distance formula}, for bounded linear functionals on a Banach space.

\begin{theorem}\label{theorem:functional-1}
Let $ \mathbb{X} $ be a reflexive Banach space and let $ f, g \in \mathbb{X}^{*} $ be linearly independent. Let $ \lambda_0 \in \mathbb{R}. $ Then the following are equivalent:\\
\item[(i)] $ (f-\lambda_0 g) \perp_B g, $
\item[(ii)] $  M_{f-\lambda_0 g} \bigcap \mathcal{N}(g) \neq \emptyset,$ 
\item[(iii)] $ \lambda_0 g $ is a best approximation to $ f $ out of $ \textit{span}\{g\}. $\\
Moreover,
$$\textit{dist}(f, \textit{span}\{ g \})=\max\{|f(x)| : x \in \mathcal{N}(g) \bigcap S_{\mathbb{X}}\} = \|f\big\rvert_{\mathcal{N}(g)}\|.$$
\end{theorem}

\begin{proof}
The equivalence of (i) and (iii) follows from the corresponding definitions, as before. To prove that each of them is equivalent to (ii), we apply Theorem $ 2.13 $ of \cite{SPM} to conclude that there exist $ z, w \in M_{f-\lambda_0 g} $ such that $ (f-\lambda_0 g)z.g(z) \geq 0  $ and  $ (f-\lambda_0 g)w.g(w) \leq 0. $ By virtue of Proposition \ref{prop:functional}, $ M_{f-\lambda_0 g} = \pm D, $ where $ D $ is a connected subset of $ S_{\mathbb{X}}. $ It is now easy to deduce that there exists $ x \in M_{f-\lambda_0 g} $ such that $ (f-\lambda_0 g)x.g(x) = 0. $ Clearly, $ f-\lambda_0 g \neq \theta. $ Therefore, we conclude that $ g(x) = 0, $ or, equivalently, $ x \in M_{f-\lambda_0 g} \bigcap \mathcal{N}(g). $ On the other hand, if $ (ii) $ holds, then taking $ x \in  M_{f-\lambda_0 g} \bigcap \mathcal{N}(g), $ we obtain that for any $ \lambda \in \mathbb{R}, $  
\[ \| (f-\lambda_0 g) + \lambda g \| \geq | (f-(\lambda_0-\lambda) g)x 
|  = | (f-\lambda_0 g)x | = \| f- \lambda_0 g\|. \]
This proves that $ (f-\lambda_0 g) \perp_B g, $ and completes the proof of the first part of the theorem. The second part follows from Theorem \ref{theorem:distance formula}. We just need to observe that since $ \mathbb{Y} = \mathbb{R}, $ the only s.i.p. on $ \mathbb{Y} $ is given by the usual multiplication of real numbers. Therefore, if $ x \in S_{\mathbb{X}} $ is such that $ g(x).y=0 $ for some $ y=\pm 1, $ then $ x \in \mathcal{N}(g). $  This establishes the theorem.
\end{proof}

\begin{remark}
Slightly digressing from our main topic of interest, we would like to note that Proposition \ref{prop:functional} allows us to improve Theorem $ 2.13 $ of \cite{SPM} by proving that the condition of strict convexity in the said theorem is redundant.
\end{remark}

Our next result on best approximations to functionals is analogous to Theorem \ref{theorem:n-dim}. The proof is omitted as it can be completed using similar arguments as before.

\begin{theorem}\label{theorem:functional-2}
Let $ \mathbb{X} $ be a reflexive Banach space and let $ n \geq 2. $ Let $ f, g_1, g_2, \ldots, g_n \in \mathbb{X}^{*} $ be such that $ g_1, g_2, \ldots, g_n $ are linearly independent and $ f \notin \mathbb{Z}, $ where $ \mathbb{Z} = \textit{span}\{ g_1, g_2, \ldots, g_n \}. $ Let $ \alpha_i \in \mathbb{R}, $ where $ 1 \leq i \leq n. $  Then $ \sum_{i=1}^{n} \alpha_i g_i $ is a best approximation to $ f $ out of $ \mathbb{Z} $ if and only if given any $ n $ scalars $ \beta_1, \beta_2, \ldots, \beta_n, $ there exist scalars $  \gamma_1, \gamma_2, \ldots, \gamma_n $
 such that the following holds true:
\begin{align*}
\| f - \sum_{i=1}^{n} \beta_i g_i \| & = \max\{ | (f-\sum_{i=1}^{n} \beta_i g_i)x | : x \in \mathcal{N}(\sum_{i=1}^{n} \gamma_i g_{i}) \bigcap  S_{\mathbb{X}}  \} \\
 & \geq \max\{ | (f-\sum_{i=1}^{n} \alpha_i g_i)x | : x \in \mathcal{N}(\sum_{i=1}^{n} \gamma_i g_{i}) \bigcap S_{\mathbb{X}} \} \\
 & = \| f - \sum_{i=1}^{n} \alpha_i g_i \|.
\end{align*}
\end{theorem}

If in addition, $ \mathbb{X} $ is strictly convex, then we have the following refinement of the above theorem: 

\begin{theorem}\label{theorem:functional-3}
Let $ \mathbb{X} $ be a reflexive strictly convex Banach space and let $ n \in \mathbb{N}. $ Let $ f, g_1, g_2, \ldots, g_n \in \mathbb{X}^{*} $ be such that $ g_1, g_2, \ldots, g_n $ are linearly independent and $ f \notin \mathbb{Z}, $ where $ \mathbb{Z} = \textit{span}\{ g_1, g_2, \ldots, g_n \}. $ Let $ \alpha_i \in \mathbb{R}, $ where $ 1 \leq i \leq n. $  Then $ \sum_{i=1}^{n} \alpha_i g_i $ is a best approximation to $ f $ out of $ \mathbb{Z} $ if and only if $ \bigcap_{i=1}^{n} \mathcal{N}(g_i) \bigcap  M_{f-\sum_{i=1}^{n} \alpha_i g_i} = \{\pm x_0\}, $ for some $ x_0 \in S_{\mathbb{X}}. $
\end{theorem}

\begin{proof}
To prove the sufficient part of the theorem, simply observe that for any $ 1 \leq j \leq n, $ and for any $ \lambda \in \mathbb{R}, $ we have the following:
\begin{align*}
\| (f - \sum_{i=1}^{n} \alpha_i g_i) + \lambda g_j \| & \geq  | \{(f-\sum_{i=1}^{n} \alpha_i g_i) + \lambda g_j \}x_0 |  \\
 & = | (f - \sum_{i=1}^{n} \alpha_i g_i)x_0 | \\
 & = \| f - \sum_{i=1}^{n} \alpha_i g_i \|.
\end{align*}

This shows that for any $ 1 \leq j \leq n, $ we have that $ (f - \sum_{i=1}^{n} \alpha_i g_i) \perp_B g_j. $ As $ \mathbb{X} $ is reflexive and strictly convex, it is easy to see that $ \mathbb{X}^* $ is smooth and consequently, Birkhoff-James orthogonality is right additive in $ \mathbb{X}^* $. Therefore, it follows that $ (f - \sum_{i=1}^{n} \alpha_i g_i) \perp_B \mathbb{Z}. $ This is clearly equivalent to the desired conclusion that $ \sum_{i=1}^{n} \alpha_i g_i $ is a best approximation to $ f $ out of $ \mathbb{Z}. $ Let us now prove the necessary part of the theorem. Clearly, for each $ 1 \leq j \leq n, $ it follows that $ (f - \sum_{i=1}^{n} \alpha_i g_i) \perp_B g_j. $ Therefore, Theorem \ref{theorem:functional-1} implies that for each $ 1 \leq j \leq n, $ $ M_{f - \sum_{i=1}^{n} \alpha_i g_i} \bigcap \mathcal{N}(g_j) \neq \emptyset. $ We also note that since $ \mathbb{X} $ is strictly convex, there exists a unique $ x_0 \in S_{\mathbb{X}} $ such that $ M_{f - \sum_{i=1}^{n} \alpha_i g_i} = \{ \pm x_0 \}. $ Combining these two observations, we obtain that $ \bigcap_{i=1}^{n} \mathcal{N}(g_i) \bigcap  M_{f-\sum_{i=1}^{n} \alpha_i g_i} = \{\pm x_0\}. $ This establishes the theorem.
\end{proof}

As an immediate application of the above theorem, we have the following elegant distance formula in the dual of a reflexive strictly convex Banach space.

\begin{theorem}\label{theorem:functional-distance formula}
Let $ \mathbb{X} $ be a reflexive strictly convex Banach space and let $ n \in \mathbb{N}. $ Let $ f, g_1, g_2, \ldots, g_n \in \mathbb{X}^{*} $ be such that $ g_1, g_2, \ldots, g_n $ are linearly independent and $ f \notin \mathbb{Z}, $ where $ \mathbb{Z} = \textit{span}\{ g_1, g_2, \ldots, g_n \}. $ Then 
$$ \textit{dist}(f, \mathbb{Z})=\max\{|f(x)| : x \in \bigcap_{i=1}^{n} \mathcal{N}(g_i) \bigcap S_{\mathbb{X}} \}. $$
\end{theorem}

\begin{proof}
A standard compactness argument shows that there exists $ \alpha_i \in \mathbb{R}, $ where $ 1 \leq i \leq n, $ such that  $ \sum_{i=1}^{n} \alpha_i g_i $ is a best approximation to $ f $ out of $ \mathbb{Z}. $ Using the necessary part of Theorem \ref{theorem:functional-3}, we obtain that $ \bigcap_{i=1}^{n} \mathcal{N}(g_i) \bigcap  M_{f-\sum_{i=1}^{n} \alpha_i g_i} = \{\pm x_0\}, $ for some $ x_0 \in S_{\mathbb{X}}. $ We also note that $ \| f-\sum_{i=1}^{n} \alpha_i g_i \| = | (f-\sum_{i=1}^{n} \alpha_i g_i)x_0 | = | f(x_0) | = \max\{|f(x)| : x \in \bigcap_{i=1}^{n} \mathcal{N}(g_i) \bigcap S_{\mathbb{X}} \}, $ where the last equality is a  consequence of the fact that $  f = f-\sum_{i=1}^{n} \alpha_i g_i, $ when restricted to $ \bigcap_{i=1}^{n} \mathcal{N}(g_i). $ Indeed, if there exists $ y_0 \in \bigcap_{i=1}^{n} \mathcal{N}(g_i) \bigcap S_{\mathbb{X}} $ such that $ | f(y_0) | > | f(x_0) | $ then we have that 
$$ \| f-\sum_{i=1}^{n} \alpha_i g_i \| \geq | (f-\sum_{i=1}^{n} \alpha_i g_i)y_0 | = | f(y_0) | > | f(x_0) | = \| f-\sum_{i=1}^{n} \alpha_i g_i \|, $$
a contradiction. As $ \textit{dist}(f, \mathbb{Z}) = \| f-\sum_{i=1}^{n} \alpha_i g_i \|, $ this completes the proof of the theorem.
\end{proof}

The above theorem, in conjunction with Theorem \ref{theorem:functional-3},  turns out to be extremely useful in studying best approximations in a  reflexive, smooth and strictly convex Banach space $ \mathbb{X}. $ To illustrate the main idea, let us begin with the following standard best approximation problem: \\

\noindent \textbf{Problem:} Let $ \mathbb{X} $ be a reflexive, smooth and strictly convex Banach space. Let $ x_0, y_1, \ldots, y_n \in \mathbb{X}, $ where $ n \geq 1, $ be such that $ y_1, \ldots, y_n $ are linearly independent and $ x_0 \notin \mathbb{Y} = \textit{span}\{y_1, \ldots, y_n\}. $ Find the (unique) best approximation to $ x_0 $ out of $ \mathbb{Y} $ and compute $ \textit{dist} (x_0, \mathbb{Y}). $\\

Theorem \ref{theorem:functional-3} and Theorem \ref{theorem:functional-distance formula} allow us to approach the above problem in the following  three steps:\\

\textit{Step 1:} We identify $ \mathbb{X} $ with its double dual $ \mathbb{X}^{**} $ via the canonical isometric isomorphism $ \psi. $ Let $ \psi(x_0) = f_0 $ and $ \psi(y_i) = g_i, $ where $ 1 \leq i \leq n. $ Let $ \mathbb{Z} = \textit{span}\{g_1, \ldots, g_n\}.  $ The original problem is clearly equivalent to finding the best approximation to $ f_0 $ out of $ \mathbb{Z} $ and computing $ \textit{dist} (f_0, \mathbb{Z}). $ As $ \mathbb{X}^{*} $ is reflexive, smooth and strictly convex, we are in a position to apply Theorem \ref{theorem:functional-3} and Theorem \ref{theorem:functional-distance formula}.\\

\textit{Step 2:} Let $ \mathcal{W} = \bigcap_{i=1}^{n} \mathcal{N}(g_i). $ Applying Theorem \ref{theorem:functional-distance formula}, we obtain that 
$$ \textit{dist} (x_0, \mathbb{Y}) = \textit{dist} (f_0, \mathbb{Z}) = \max\{|f_0(x)| : x \in \mathcal{W} \bigcap S_{\mathbb{X}^{*}} \}. $$ 

\textit{Step 3:} Since $ \mathbb{X}^{*} $ is strictly convex, so is $ \mathcal{W}. $ Therefore, $ f_0 \big\rvert_{\mathcal{W}} $ attains norm at only one pair of points, say, $ \pm  h_0 \in S_{\mathcal{W}}. $  Applying Theorem \ref{theorem:functional-3}, $ \sum_{i=1}^{n} \alpha_i g_i $ is the unique best approximation to $ f_0 $ out of $ \mathbb{Z} $ if and only if $ M_{f_0-\sum_{i=1}^{n} \alpha_i g_i} = \{ \pm h_0 \}. $  We also note that since $ \mathbb{X}^{*} $ is smooth, there exists a unique $ \xi_0 \in \mathbb{X}^{**} $ such that $ M_{\xi_0} = \{ \pm h_0 \}. $ Therefore, we must have that $ f_0-\sum_{i=1}^{n} \alpha_i g_i = \xi_0, $ which completely describes the unique best approximation to $ f_0 $ out of $ \mathbb{Z} $ by means of the following equation:
$$ \sum_{i=1}^{n} \alpha_i g_i = f_0 - \xi_0. $$
We have therefore obtained a complete answer to the original problem in light of the fact that $ \sum_{i=1}^{n} \alpha_i y_i $ is the unique best approximation to $ x_0 $ out of $ \mathbb{Y}. $ \\

\begin{remark}
Unlike Theorem \ref{theorem:functional-3} and Theorem \ref{theorem:functional-distance formula}, we have made use of the strict convexity of $ \mathbb{X}^* $ in the above algorithm. Since $ \mathbb{X} $ is reflexive, this is ensured by (and is in fact equivalent to) the smoothness of $ \mathbb{X}. $
\end{remark}

The above algorithm can be applied in a more efficient manner from computational point of view, in case of $ \ell_{p}^{n} $ spaces, where $ 1 < p < \infty. $ This is because of the well-known identification of $ \ell_{p}^{n} $ with $ \ell_{q}^{{n}^{*}}, $ where $ \frac{1}{p} + \frac{1}{q} = 1, $ under the mapping $ \Omega : \ell_p^n \longrightarrow \ell_{q}^{{n}^{*}} $ given by
\[ \Omega (\widetilde{x}) = f \qquad \forall ~\widetilde{x} =  (x_1, x_2, \ldots, x_n) \in \ell_p^n, \]
where $ f : \ell_q^n \longrightarrow \mathbb{R} $ is given by 
\[ f(e_i) = x_i ~ \forall~ i \in \mathbb{N}, \] 
$ \{ e_i : 1 \leq i \leq n \} $ being the standard ordered basis for $ \ell_{q}^n. $ Let us explain this in more detail in the following:\\

\noindent \textbf{Problem:} Let $ \widetilde{x}, \widetilde{y_1}, \ldots, \widetilde{y_m} \in \ell_p^n, $ where $ 1 < p < \infty $ and $ 1 \leq m < n, $ be such that $ \widetilde{y_1}, \ldots, \widetilde{y_m} $ are linearly independent and $ \widetilde{x} \notin \mathbb{Y} = \textit{span}\{\widetilde{y_1}, \ldots, \widetilde{y_m
}\}. $ Compute $ \textit{dist} (\widetilde{x}, \mathbb{Y}). $\\

\textit{Step 1:} Let $ \widetilde{x}=(x_1, \ldots, x_n) $ and let $ \widetilde{y_i}= (y_{i1}, \ldots, y_{in}), $ for each $ 1 \leq i \leq m. $ Let $ \Omega(\widetilde{x}) = f $ and let $ \Omega(\widetilde{y_i}) = g_i,  $ for each $ 1 \leq i \leq m. $ Let $ \mathbb{Z} = \textit{span}\{g_1, \ldots, g_m\}.  $ The original problem is clearly equivalent to finding the best approximation to $ f $ out of $ \mathbb{Z} $ and computing $ \textit{dist} (f, \mathbb{Z}). $\\

\textit{Step 2:} Let $ \mathcal{W} = \bigcap_{i=1}^{m} \mathcal{N}(g_i). $ Clearly, $ \mathcal{W} $ is a subspace in $ \mathbb{R}^n $ determined by the fact that $ (\eta_1, \ldots, \eta_n) \in \mathcal{W} $ if and only if $ (\eta_1, \ldots, \eta_n) $ satisfies the following system of equations:
\[ \sum_{j=1}^{n} y_{kj}\eta_j = 0~ ;~ 1 \leq k \leq m. \]

Applying Theorem \ref{theorem:functional-distance formula}, we obtain that 
$$ \textit{dist} (\widetilde{x}, \mathbb{Y}) = \textit{dist} (f, \mathbb{Z}) = \max\{|f(x)| : x \in \mathcal{W} \bigcap S_{\ell_{q}^{n}} \}. $$\\

In particular, when $ n=2, $ we have the following explicit distance formula:

\begin{theorem}\label{theorem:distance formula2} Let $ \mathbb{X} = \ell_p^2, $ where $ 1 < p < \infty. $ Let $ \widetilde{x} = (a, b) \in \mathbb{X} $, let $(c,d)\neq (0,0)$ and let $ \mathbb{Y} = \textit{span}\{ (c, d) \} $ be a one-dimensional subspace in $ \mathbb{X} $ be such that $ \widetilde{x} \notin \mathbb{Y}. $ Then 
\begin{align*}
\textit{dist} (\widetilde{x}, \mathbb{Y}) = \frac{|ad-bc|}{(|c|^q+|d|^q)^{\frac{1}{q}}},
\end{align*}
where $ \frac{1}{p} + \frac{1}{q} = 1. $
\end{theorem}

\begin{proof}
Let $ \widetilde{y} = (c, d). $ Let $ \Omega : \ell_p^2 \longrightarrow \ell_{q}^{{2}^{*}} $ be the canonical isometric isomorphism. Let $ \Omega(\widetilde{x}) = f $ and let $ \Omega(\widetilde{y}) = g. $ Let us first assume that $ d \neq 0. $ Now, we have that 
\[ \mathcal{N}(g) = \{(\eta_1, \eta_2) \in \mathbb{R}^2 : c\eta_1+d\eta_2=0 \}.\]
Since $ d \neq 0, $ it follows that $ (1, -\frac{c}{d}) \in \mathcal{N}(g). $ Let $ \alpha \in \mathbb{R} $ be such that $ \| \alpha (1, -\frac{c}{d}) \|_{q} = 1. $ An easy calculation implies that $ | \alpha | = \frac{1}{(1+|\frac{c}{d}|^q)^{\frac{1}{q}}}. $ Applying Theorem \ref{theorem:functional-distance formula}, we obtain that 
\[ \textit{dist} (\widetilde{x}, \mathbb{Y}) = | f(\alpha (1, -\frac{c}{d})) | = \frac{|ad-bc|}{|d|(1+|\frac{c}{d}|^q)^{\frac{1}{q}}} = \frac{|ad-bc|}{(|c|^q+|d|^q)^{\frac{1}{q}}}, \]
proving the first equality. Similarly, if $ c \neq 0, $ it can be proved that 

\[  \textit{dist} (\widetilde{x}, \mathbb{Y}) = \frac{|ad-bc|}{(|c|^q+|d|^q)^{\frac{1}{q}}}. \]
This establishes the theorem.
\end{proof}

The methods developed in the present article point to the fact that every best approximation problem in a reflexive smooth and strictly convex Banach space corresponds to a maximization problem in a particular subspace of the dual space. Moreover, in case of $ \ell_p^n $ spaces, where $ 1 < p < \infty, $ computation of the distance of a given point from a given subspace reduces to a trivial calculation when the subspace is of codimension one in $ \ell_p^n. $ Indeed, in this case we just need to evaluate a particular functional on the unit sphere of a one-dimensional subspace of $ \ell_q^n, $ where $ \frac{1}{p} + \frac{1}{q} = 1. $ This is illustrated in the following example:  

\begin{example}\label{example:1}
Let $ \mathbb{X} = \ell_p^3, $ where $ 1 < p < \infty. $ Let $ x_0 = (a, b, c) $ and let $ \mathbb{Y}_0 = \textit{span}\{ (1, 0 , -1), (1, 2, 1) \}, $ where $ a, b, c \in \mathbb{R}. $ Let us compute $ \textit{dist} (x_0, \mathbb{Y}_0) $ by using the ideas developed above. Let $ \Omega  : \ell_p^3 \longrightarrow \ell_{q}^{{3}^{*}} $ be the canonical isometric isomorphism, where $ \frac{1}{p} + \frac{1}{q} = 1. $ Let $ y_1 = (1, 0 , -1) $ and let $ y_2 = (1, 2 , 1). $ Let us assume that $ \Omega(x_0)=f_0,~  \Omega(y_1)=g_1,~\Omega(y_2)=g_2. $ An easy computation shows that $ \mathcal{N}(g_1) \bigcap \mathcal{N}(g_2) = \textit{span}\{ (1, -1, 1) \}. $ Since $ \|(1, -1, 1)\|_q = 3^{\frac{1}{q}}, $ applying Theorem \ref{theorem:functional-distance formula} we obtain that 
\[ \textit{dist} (x_0, \mathbb{Y}_0) = | f_0 (3^{\frac{-1}{q}}(1, -1, 1)) | = 3^{\frac{-1}{q}} | a-b+c |.  \]   
\end{example}

Let us end this section with the following remark, emphasizing the utility of the methods developed here in studying best approximation problems.

\begin{remark}
The above example, when viewed purely in the $ \ell_p^3 $ setting, amounts to solving the following minimization problem:
\[ \min_{\alpha, \beta \in \mathbb{R}} \{ | a-\alpha-\beta |^p + | b-2\beta |^p + | c+\alpha-\beta |^p \}^{\frac{1}{p}}, \]
which require some computational efforts. However, applying the duality theory of best approximations developed in this article, we can readily conclude via a trivial computation that the answer to the above minimization problem is given by $ 3^{\frac{-1}{q}} | a-b+c |, $ where $ \frac{1}{p} + \frac{1}{q} = 1. $ Equivalently, this can be expressed by means of the following inequality:
\[ | a-\alpha-\beta |^p + | b-2\beta |^p + | c+\alpha-\beta |^p \geq 3^{1-p} | a-b+c |^p ~\forall~ \alpha, \beta, a, b, c \in \mathbb{R}. \]
Indeed, the methods developed in this article in view of best approximations in $ \ell_p $ spaces give rise to a family of such inequalities, which seem not so easy to prove otherwise. It should be noted that Theorem \ref{theorem:functional-distance formula} guarantees that each of these inequalities is optimal.
\end{remark}

\section{A comparison with the classical duality principle} The aim of this section is to make a comparative analysis of some of our results, namely, Theorem \ref{theorem:functional-3} and Theorem \ref{theorem:functional-distance formula}, with the following well-known classical duality principle which is one of the very basic and foundational results in approximation theory:

\[ \textit{dist} (x, \mathbb{Y}) = \inf \{ \| x - y \|_{\mathbb{X}} : y \in \mathbb{Y} \} = \sup\{ f(x) : \| f \|_{\mathbb{X}^*} \leq 1, f\lvert_{\mathbb{Y}} = 0\}, \]

where $ \mathbb{Y} $ is a subspace of $ \mathbb{X}. $\\

We begin with the observation that neither Theorem \ref{theorem:functional-3} nor Theorem \ref{theorem:functional-distance formula} of this article follows as a direct consequence of the above duality principle. In particular, the use of Birkhoff-James orthogonality in studying best approximations in a dual space allows us to approach the problem from a point of view which is computationally much more convenient. To explain this further, let us first have a closer look at the above distance formula given by the classical duality principle. To explicitly compute the distance between a given point $ x $ and a given subspace $ \mathbb{Y} $ of the Banach space $ \mathbb{X} $ by applying the above formula, the major hindrance is that in general we do not have an explicit description of all those functionals $ f $ in $ B_{\mathbb{X}^*} $ which vanish on $ \mathbb{Y}. $ This can be remedied from a theoretical point of view in case $ \mathbb{X} $ is reflexive and $ \mathbb{Y} $ is of codimension one in $ \mathbb{X}. $ Since $ \mathbb{X} $ is reflexive, it follows from Theorem $ 1 $ of \cite{Jb} that there exists $ z \in S_{\mathbb{X}} $ such that $ z \perp_B \mathbb{Y}. $ It is now easy to see that the following holds true:
\[ \{ f \in S_{\mathbb{X}^*} : f(\mathbb{Y}) = 0 \} = \{ f \in S_{\mathbb{X}^*} : |f(z)| = 1 \}. \]

If in addition, $ \mathbb{X} $ is smooth then there exists a unique $ f_0 \in S_{\mathbb{X}^*} $ such that $ f_0(z)=1. $ Therefore, by applying the classical duality principle, we obtain that
\[ \textit{dist} (x, \mathbb{Y}) = | f_0(x) |. \]

Unfortunately, Theorem $ 1 $ of \cite{Jb} does not provide us with a method to compute $ z $ such that $ z \perp_B \mathbb{Y}. $ In fact, even in the finite-dimensional case, it is computationally not straightforward to explicitly find out such a $ z $ in a general Banach space, when the norm is not induced by an inner product. The reader is invited to verify this claim for $ \ell_p^n $ spaces, where $ n \geq 2 $ and $ p \in (1, \infty) \setminus \{ 2 \}. $ Indeed, to find out such a $ z $ such that $ z \perp_B \mathbb{Y}, $ it is readily seen from the expression of the unique s.i.p. on $ \ell_p^n $ that the problem amounts to solving \emph{nonlinear} equations, which is quite well-known to be computationally hard. Moreover, if the codimension of $ \mathbb{Y} $ in $ \mathbb{X} $ is strictly greater than one, then it is easy to see that $ \mathbb{Y} $ cannot be written as $ \mathbb{Y} = v^{\perp}, $ for any $ v \in S_{\mathbb{X}}. $ Therefore, it becomes an even more difficult problem to have an explicit description of all those functionals $ f $ in $ B_{\mathbb{X}^*} $ which vanish on $ \mathbb{Y}. $ This explains the computational difficulty in applying the classical duality principle in determining the distance between a given point $ x $ and a given subspace $ \mathbb{Y} $ of the Banach space $ \mathbb{X}. $ \\

In stark contrast to the above, in case of a reflexive and strictly convex Banach space $ \mathbb{X} $, our method allows us to approach the problem in a computationally efficient manner, by identifying $ \mathbb{X} $ with $ \mathbb{X}^{**} $ and then by applying Theorem \ref{theorem:functional-distance formula}. Indeed, in case $ \mathbb{X} $ is a finite-dimensional strictly convex Banach space, and $ \mathbb{Y} $ is of codimension one in $ \mathbb{X}, $ the problem of computing $ \textit{dist} (x, \mathbb{Y}) $ reduces to a trivial computation involving the solution to a given system of \emph{linear} equations. This explains the computational advantage obtained in applying Theorem \ref{theorem:functional-distance formula} over the classical duality principle and it can be further verified readily by trying to establish Theorem \ref{theorem:distance formula2} or Example \ref{example:1}, by the later alone.

\begin{acknowledgement}
The author feels elated to acknowledge the delightful company of his beloved childhood friend Dr. Chandan Das, a humanitarian physician with an empathetic approach. The author would like to thank the referees for the thoughtful comments and suggestions which led to definite improvement of the paper.
\end{acknowledgement}

\bibliographystyle{amsplain}

\end{document}